\documentclass[a4paper,11pt,reqno]{amsproc}
\usepackage{amssymb}
\usepackage[width=132mm,height=190mm]{geometry} 
\usepackage{enumerate}
\usepackage{lmodern}
\newtheorem{corollary}{Corollary}
\newtheorem{theorem}[corollary]{Theorem}
\newtheorem{fact}[corollary]{Fact}
\newtheorem{proposition}[corollary]{Proposition}
\newtheorem{lemma}[corollary]{Lemma}
\theoremstyle{definition}
\newtheorem{definition}[corollary]{Definition}
\newtheorem{remark}[corollary]{Remark}

\newcommand{\vek}[1]{\mathbf{#1}} 
\newcommand{\mat}[1]{\mathbf{#1}}
\newcommand{\abs}[1]{\lvert#1\rvert}
\newcommand{\whom}{\mathrm{w}_{\mathrm{hom}}}

\newcommand{\wham}{\mathrm{w}_{\mathrm{Ham}}}
\newcommand{\wlee}{\mathrm{w}_{\mathrm{Lee}}}

\DeclareMathOperator{\Hom}{Hom}
\DeclareMathOperator{\PG}{PG}
\DeclareMathOperator{\GL}{GL}

\newcommand{\points}{\mathcal{P}}

\newcommand{\Z}{\mathbb{Z}}
\newcommand{\N}{\mathbb{N}} 
\newcommand{\Q}{\mathbb{Q}}
\newcommand{\C}{\mathbb{C}}
\newcommand{\F}{\mathbb{F}}

\newcommand{\wrt}{w.\,r.\,t.}  
\newcommand{\ie}{i.\,e.}

\title{Further Results on Homogeneous Two-Weight Codes}

\author{Thomas Honold}
\address{Thomas Honold\\
  Institute of Information and Communication Engineering\\
  Zhejiang University, Zheda Road\\
  310027 Hangzhou, China}
\email{honold@zju.edu.cn}
\thanks{Reprint of the conference paper published in the 
  Proceedings of the Fifth International Workshop on Optimal
  Codes and Related Topics (OC2007), White Lagoon, Bulgaria, June
2007, pp.\ 80--86.}

\date{}

\begin{document}

\begin{abstract}
  The results of \cite{emt:oc4graphs,emt:graphs} on linear homogeneous
  two-weight codes over finite Frobenius rings are exended in two
  ways: It is shown that certain non-projective two-weight codes give
  rise to strongly regular graphs in the way described in
  \cite{emt:oc4graphs,emt:graphs}. Secondly, these codes are used to
  define a dual two-weight code and strongly regular graph similar to
  the classical case of projective linear two-weight codes over finite
  fields \cite{calderbank-kantor86}.
\end{abstract}



\maketitle

\section{Introduction}\label{sec:intro}

A finite ring $R$ is said to be a Frobenius ring if there exists a character
$\chi\in\widehat{R}=\Hom_\Z(R,\C^\times)$ whose kernel contains no
nonzero left (or right) ideal of $R$.
The (normalized) homogeneous weight $\whom\colon R\to\C$ on a finite
Frobenius ring $R$ is defined by
\begin{equation}
  \label{eq:whom}
  \whom(x)=1-\frac{1}{\abs{R^\times}}\sum_{u\in R^\times}\chi(ux).
\end{equation}
(This does not depend on the choice of $\chi$.)
The function $\whom$ is the unique complex-valued function on $R$ satisfying
$\whom(0)=0$, $\whom(ux)=\whom(x)$ for $x\in R$, $u\in R^\times$ and
$\sum_{x\in I}\whom(x)=\abs{I}$ for all nonzero left ideals
$I\leq{}_RR$ (and their right counterparts).

The homogeneous weight on a finite Frobenius ring is a generalization of both
the Hamming weight on $\F_q$ ($\whom(x)=\frac{q}{q-1}\wham(x)$ for
$x\in\F_q$) and the Lee weight on $\Z_4$ ($\whom(x)=\wlee(x)$ for
$x\in\Z_4$). It was introduced in
\cite{ioana-werner97} for the case $R=\Z_m$ and generalized to
Frobenius rings in \cite{greferath-schmidt99,st:homippi}.

In \cite{emt:oc4graphs,emt:graphs} it was shown that a linear code $C$
over a finite Frobenius ring with exactly two nonzero homogeneous
weights and satisfying certain nondegeneracy conditions gives rise to
a strongly regular graph with $C$ as its set of vertices. In the
classical case $R=\F_q$ this result has been known for a long time and
forms part of a more general correspondence between projective
linear $[n,k]$ two-weight codes over $\F_q$ and certain strongly regular Cayley
graphs of $(\F_q^k,+)$ resp.\ regular partial difference sets
in $(\F_q^k,+)$, and their (appropriately defined) duals
(cf.\cite{calderbank-kantor86,delsarte72}).

The purpose of this work is to generalize the results of
\cite{emt:oc4graphs,emt:graphs} to a larger class of homogeneous
two-weight codes (so-called modular two-weight codes) and establish
for these codes the classical correspondence (Theorems~3.2 and~5.7 of
\cite{calderbank-kantor86}) in full generality. 

\section{A Few Properties of Frobenius Rings and their Homogeneous
  Weights}\label{sec:properties}

For a subset $S$ of a ring $R$
let ${}^\perp S=\{x\in R;xS=0\}$, $S^\perp=\{x\in R;Sx=0\}$. Similarly,
for $S\subseteq R^n$ let ${}^\perp S=\{\vek{x}\in R^n;\vek{x}\cdot S=0\}$ and
$S^\perp=\{\vek{x}\in R^n;S\cdot\vek{x}=0\}$, where
$\vek{x}\cdot\vek{y}=x_1y_1+\dots+x_ny_n$.

\begin{proposition}
  \label{prop:same-shape}
  A finite ring $R$ is a Frobenius ring iff for every matrix
  $\mat{A}\in R^{m\times n}$ the left row space
  $C=\{\vek{x}\mat{A};\vek{x}\in R^m\}$ and the right column space
  $D=\{\mat{A}\vek{y};\vek{y}\in R^n\}$ have the same cardinality.
\end{proposition}

From now on we suppose that $R$ is a finite Frobenius ring with
homogeneous weight $\whom$. 

First we determine the set of all $x\in R$ satisfying $\whom(x)=0$.
Let $S_i=Rs_i$, $1\leq i\leq\tau$, be the different left ideals of $R$
of order $2$ and $S=S_1+\dots+S_\tau$.  The set $S$ is a two-sided
ideal of $R$ of order $2^\tau$, whose elements are the subset sums of
$\{s_1,\dots,s_\tau\}$.
Define $S_0\subseteq S$ as the set of all sums of an even number of
elements from $\{s_1,\dots,s_\tau\}$ (``even-weight subcode of $S$'').
Note that $S_0$ is a subgroup of $(R,+)$, trivial for
$\tau\leq 1$ and nontrivial (of order $2^{\tau-1}$) for $\tau\geq 2$.
\begin{proposition}
  \label{prop:S_0}
  We have $\whom(x)\geq 0$ for all $x\in R$, and $\{x\in
  R;\whom(x)=0\}=S_0$. Moreover, $\whom(x+y)=\whom(x)$ for all $x\in R$
  and $y\in S_0$.
\end{proposition}
\begin{fact}[{\cite[Th.~2]{wt:egal}}]
  \label{fact:egal}
\begin{equation}
  \label{eq:egal}
  \sum_{x\in I}\whom(x+c)=\abs{I}
\end{equation}
for all nonzero left (or right) ideals $I$ of $R$ and all $c\in
R$. 
\end{fact}
The following correlation property of $\whom$ turns out to be crucial.
\begin{proposition}
  \label{prop:whom-corr}
  For a nonzero left ideal $I$ of $R$ and $r,s\in R$ we have
  \begin{equation}
    \label{eq:whom-corr}
    \sum_{x\in I}\whom(x)\whom(xr+s)
    =
    \begin{cases}
\abs{I}+\abs{I}\cdot\frac{\abs{R^\times\cap(1+I^\perp)}}
    {\abs{R^\times}}\cdot\bigl(1-\whom(s)\bigr)&\text{if $\abs{Ir}=\abs{I}$},\\
    \abs{I}&\text{if $\abs{Ir}<\abs{I}$}.
    \end{cases}
 \end{equation}
  In particular $\sum_{x\in
    R}\whom(x)^2=\abs{R}+\frac{\abs{R}}{\abs{R^\times}}$.
\end{proposition}
For vectors $\vek{x},\vek{y}\in R^k$ we write $\vek{x}\sim\vek{y}$ if
$\vek{x}R^\times=\vek{y}R^\times$. By \cite[Prop.~5.1]{wood99a}
this is equivalent to $\vek{x}R=\vek{y}R$.
\begin{proposition}
  \label{prop:whom-corr-vek}
  For nonzero words $\vek{g},\vek{h}\in R^k$ and $s\in R$ we have
  \begin{equation}
    \label{eq:whom-corr-vek}
    \sum_{\vek{x}\in R^k}\whom(\vek{x}\cdot\vek{g})
    \whom(\vek{x}\cdot\vek{h}+s)=
    \begin{cases}
      \abs{R}^k+\frac{\abs{R}^k}{\abs{\vek{g}R^\times}}
      \cdot\bigl(1-\whom(s)\bigr) 
      &\text{if $\vek{g}\sim\vek{h}$},\\
      \abs{R}^k
      &\text{if $\vek{g}\nsim\vek{h}$}.
    \end{cases}
  \end{equation}
\end{proposition}

\section{Modular Two-Weight Codes, Partial Difference Sets and Strongly
  Regular Cayley Graphs}

Given a positive integer $k$, the set of nonzero cyclic submodules of the free
right module $R^k_R$ is denoted by $\points$. The elements of
$\points$ are referred to as \emph{points} of the projective geometry
$\PG(R^k_R)$, and a multiset $\alpha\colon\points\to\N_0$ is referred
to as a \emph{multiset in $\PG(R^k_R)$}.

With a left linear code $C\leq{}_RR^n$ generated by $k$ (or fewer)
codewords and having no all-zero coordinate we associate a multiset
$\alpha_C$ in $\PG(R^k_R)$ of cardinality $n$ in the following way: If
$C=\{\vek{x}\mat{G};\vek{x}\in R^k\}$ with
$\mat{G}=(\vek{g}_1|\vek{g}_2|\dots|\vek{g}_n)\in R^{k\times n}$,
define $\alpha_C\colon\points\to\N_0$ by
$\alpha(\vek{g}R)=\abs{\{j;1\leq j\leq n\wedge
  \vek{g}_jR=\vek{g}R\}}$. The relation $C\leftrightarrow\alpha_C$
defines a bijection between classes of monomially isomorphic left
linear codes over $R$ generated by $k$ codewords and orbits of the
group $\GL(R^k_R)$ on multisets in $\PG(R^k_R)$.

\begin{definition}
  \label{dfn:modular}
  
  A linear code $C\leq{}_RR^n$ is said to be \emph{modular} if there exists
  $r\in\Q$ such that for all points $\vek{g}R$ of $\PG(R^k_R)$ either
  $\alpha_C(\vek{g}R)=0$ or
  $\alpha_C(\vek{g}R)=r\abs{\vek{g}R^\times}$. The number $r$ is
  called the \emph{index} of $C$.
\end{definition}

The property of $C$ described in Def.~\ref{dfn:modular} does not
depend on the choice of $\alpha_C$ (not even on the dimension $k$).
Hence modularity of a linear code is a well-defined concept.

If $A\subseteq R^k\setminus\{\vek{0}\}$ satisfies $AR^\times=A$, the
matrix $\mat{G}$ with the vectors of $A$ as columns generates a
modular (left) linear code of length $\abs{A}$ and index $1$.


Note that projective codes over $\F_q$ are
modular of index $\frac{1}{q-1}$ and regular
projective codes over $R$ as defined in \cite{emt:oc4graphs,emt:graphs} are
modular of index $\frac{1}{\abs{R^\times}}$.

\begin{fact}[{\cite[Th.~5.4]{wood02}}]
  \label{fact:equidistant}
  A linear code $C\leq{}_RR^n$ is a one-weight code (\ie\ equidistant
  \wrt\ $\whom$)
  iff $C$ is modular and
  $\bigl\{\vek{g}\in R^k\setminus\{\vek{0}\};\alpha_C(\vek{g}R)>0\bigr\}$
  is the set of nonzero vectors of a submodule of $R^k_R$.
\end{fact}

The main purpose of this paper is a combinatorial characterization
of linear homogeneous two-weight codes over $R$, \ie\ linear codes over $R$
having exactly two nonzero homogeneous weights $w_1<w_2$. Assuming that $C$ is
such a code, we set $w_0=0$, $C_i=\{\vek{c}\in C;\whom(\vek{c})=w_i\}$ and
$b_i=\abs{C_i}$ for $i=0,1,2$. 

By Prop.~\ref{prop:S_0} we have $\whom(\vek{c})=0$ iff $\whom(c_j)=0$ for
$1\leq j\leq n$, the set $C_0$ is a subgroup of $(C,+)$ and $C_1$, $C_2$ are
  unions of cosets of $C_0$. If the weights $w_1$, $w_2$ and
  $b_0=\abs{C_0}$ are known, the frequencies $b_1$, $b_2$ can be
  computed from the equations $b_1+b_2=\abs{C}-\abs{C_0}$,
  $b_1w_1+b_2w_2=\sum_{\vek{c}\in C}\whom(\vek{c})=n\abs{C}$ (assuming
  that $C$ has no all-zero coordinate) and are given by
  \begin{equation}
    \label{eq:b_i}
    b_1=\frac{(w_2-n)\abs{C}-w_2\abs{C_0}}{w_2-w_1},\quad
    b_2=\frac{(n-w_1)\abs{C}+w_1\abs{C_0}}{w_2-w_1}.
  \end{equation}
\begin{lemma}
  \label{lma:whom-corr-code}
  For a modular code $C\leq{}_RR^n$ of index $r$ and $\vek{d}\in R^n$ we have
  \begin{equation}
    \label{eq:whom-corr-code}
    \sum_{\vek{c}\in C}\whom(\vek{c})\whom(\vek{c}+\vek{d})
    =\abs{C}\cdot\bigl(n^2+rn-r\cdot\whom(\vek{d})\bigr).
  \end{equation}
\end{lemma}
In the special case $\vek{d}=\vek{0}$ Lemma~\ref{lma:whom-corr-code}
reduces to
$\sum_{\vek{c}\in C}\whom(\vek{c})^2=(n^2+rn)\abs{C}$.
\begin{lemma}
  \label{lma:w_1w_2}
   The nonzero weights $w_1,w_2$ of a modular two-weight code
   $C\leq{}_RR^n$ of index $r$ satisfy the relation
   \begin{equation}
     \label{eq:w_1w_2}
     (w_1+w_2)n\abs{C}=(n^2+rn)\abs{C}+w_1w_2\bigl(\abs{C}-\abs{C_0}\bigr).
   \end{equation}
\end{lemma}
\begin{lemma}
  \label{lma:whom-coset-code}
  For a modular two-weight code $C\leq{}_RR^n$ of index $r$ and
  $\vek{d}\in R^n$ we have
  \begin{equation}
    \label{eq:whom-coset-code}
    \sum_{\vek{c}\in C_1}\whom(\vek{c}+\vek{d})=
    b_1w_1+\left(b_1-\frac{b_1w_1}{n}\right)\whom(\vek{d})
  \end{equation}
\end{lemma}
\begin{remark}
  \label{rmk:whom-code}
  Lemmas~\ref{lma:whom-corr-code} and ~\ref{lma:whom-coset-code} can
  be generalized to
  \begin{align*}
    \sum_{\vek{c}\in C}\whom(\vek{c})\whom(c_j+d_j)
    &=\abs{C}\cdot\bigl(n+r-r\cdot\whom(d_j)\bigr)\quad\text{and}\\
    \sum_{\vek{c}\in C_1}\whom(c_j+d_j)&=
    \frac{b_1w_1}{n}+\left(b_1-\frac{b_1w_1}{n}\right)\whom(d_j)
  \end{align*}
  respectively, where $j$ is any coordinate of $R^n$ and $d_j\in R$.
  In particular $\sum_{\vek{c}\in
    C_1}\whom(c_j)=\frac{b_1w_1}{n}$ is independent
  of $j$.
\end{remark}
Recall that a (simple)
graph $\Gamma$ is \emph{strongly regular with parameters
$(N,K,\lambda,\mu)$} if $\Gamma$ has $N$ vertices, is regular of degree $K$
and any two adjacent (resp.\ nonadjacent) vertices have $\lambda$
(resp.\ $\mu$) common neighbours. The graph $\Gamma$ is called \emph{trivial}
if $\Gamma$ or its complement is a disjoint union of cliques of the same
size. This is equivalent to $\mu=0$ resp.\ $\mu=K$.

A subset $D\subset G$ of an (additively written) abelian group $G$ is
said to be a \emph{regular $(N,K,\lambda,\mu)$ partial difference set
  in $G$} if $N=\abs{G}$, $K=\abs{D}$, $0\notin D$, $-D=D$, and the
multiset $D-D$ represents each element of $D$ exactly $\lambda$ times and
each element of $G\setminus\bigl(D\cup\{0\}\bigr)$ exactly $\mu$ times;
cf.~\cite{ma94}.

If $D$ is a regular $(N,K,\lambda,\mu)$ partial difference set
in $G$, then the graph $\Gamma(G,D)$ with vertex set $G$ and edge set
$\bigl\{\{x,x+d\};x\in G,d\in D\bigr\}$, the so-called \emph{Cayley graph} of
$G$ \wrt\ $D$, is strongly regular with parameters
$(N,K,\lambda,\mu)$.

We are now ready to generalize the main result of
\cite{emt:graphs,emt:oc4graphs} to modular two-weight codes. For a
two-weight code $C$ we denote the Cayley graph $\Gamma(C/C_0,C_1/C_0)$
by $\Gamma(C)$. Thus the vertices of $\Gamma(C)$ are the cosets of $C_0$
in $C$, and two cosets $\vek{c}+C_0$, $\vek{d}+C_0$ are
adjacent iff $\whom(\vek{c}-\vek{d})=w_1$. As we have already mentioned,
Prop.~\ref{prop:S_0} ensures that $\Gamma(C)$ is well-defined.

\begin{theorem}
  \label{thm:Gamma(C)}
  The graph $\Gamma(C)$ associated with a modular two-weight code over
  a finite Frobenius ring $R$ is strongly regular with parameters
  \begin{gather*}
    N=\frac{\abs{C}}{\abs{C_0}},\quad
    K=\frac{(w_2-n)N-w_2}{w_2-w_1},\\
    \lambda=\frac{K\left(\frac{w_1^2}{n}-2w_1\right)+w_2(K-1)}{w_2-w_1},\quad
    \mu=\frac{K\left(\frac{w_1w_2}{n}-w_1-w_2\right)+w_2K}{w_2-w_1}.
  \end{gather*}
  The graph $\Gamma(C)$ is trivial iff $w_1=n$.
\end{theorem}
\begin{remark}
  \label{rmk:trivial}
  Since $\Gamma(C)$ is a Cayley graph,
  the preceding argument shows that $\Gamma(C)$ is trivial iff the
  codewords of weight $0$ and $w_2$ form a linear subcode of $C$ (and
  the cocliques of $\Gamma(C)$ are the cosets of $(C_0+C_2)/C_0$ in
  this case).
\end{remark}

\section{The Dual of a Modular Two-Weight Code}\label{sec:dual}

Suppose $C\leq{}_RR^n$ is a two-weight code over a finite Frobenius
ring with nonzero
weights $w_1<w_2$ and frequencies $b_1$, $b_2$. Let $\mat{M}_i\in
R^{b_i\times n}$ ($i=1,2$) be matrices whose rows are the 
codewords of $C$ of weight $w_i$ in some order.
\begin{definition}
  \label{dfn:dual}
  The right linear code $C'\leq R^{b_1}_R$ generated by the columns of
  $\mat{M}_1$ is called the \emph{dual of the two-weight code $C$}.
\end{definition}
The code $C'$ is modular of index $1$ (no matter whether $C$ is modular or
not).

\begin{theorem}
  \label{thm:dual}
  If $C\leq{}_RR^n$ is a modular two-weight code with $C_0=\{\vek{0}\}$, its
  dual $C'$ is also a (modular) two-weight code with
  $C'_0=\{\vek{0}\}$ and nonzero weights
  \begin{equation}
    \label{eq:dual}
    w_1'=\frac{(w_2-n-r)\abs{C}}{w_2-w_1}=\frac{b_1w_1}{n},\quad 
    w_2'=\frac{(w_2-n)\abs{C}}{w_2-w_1}.
  \end{equation}
\end{theorem}
\begin{theorem}
  \label{thm:Gamma(C')}
  Under the assumptions of Th.~\ref{thm:dual}, the graph $\Gamma(C')$
  is strongly regular with parameters
  \begin{equation*}
    N'=\abs{C},\quad K'=\frac{n}{r},\quad
    \lambda'=\frac{2n-w_1-w_2}{r}+\frac{w_1w_2}{r^2\abs{C}},\quad
    \mu'=\frac{w_1w_2}{r^2\abs{C}}.
  \end{equation*}
  The graph $\Gamma(C')$ is trivial iff $w_1=n$ (\ie\ iff $\Gamma(C)$
  is trivial).
\end{theorem}
\begin{theorem}
  \label{thm:equivalence}
  Let $C\leq{}_RR^n$ be a modular linear code over a finite Frobenius
  ring $R$ generated by $\mat{G} =(\vek{g}_1|\dots|\vek{g}_n)\in
  R^{k\times n}$. Let 
  $D\leq R^k_R$ be the right column space of $\mat{G}$.  Suppose $C$
  has no all-zero coordinate and satisfies $C_0=\{\vek{0}\}$. Then the
  following are equivalent:
  \begin{enumerate}[(i)]
  \item $C$ is a homogeneous two-weight code;
  \item $\Omega=\vek{g}_1R^\times\cup\dots\cup\vek{g}_n R^\times$ is a
    regular partial difference set in $(D,+)$ and
    $\Omega\cup\{\vek{0}\}$ is not a submodule of $R^k_R$.
  \end{enumerate}
\end{theorem}
\begin{remark}
  \label{rmk:equivalence}
  Under the assumptions of Th.~\ref{thm:equivalence} the set
  $\Omega\cup\{\vek{0}\}$ is a submodule of $R^k_R$ iff $C$ is a
  homogeneous one-weight code, and $D\setminus\Omega$ is a submodule
  of $R^k_R$ iff $C$ is a homogeneous two-weight code with
  $w_1=n$. 
\end{remark}


\def\cprime{$'$}

\end{document}